\documentclass[preprint]{article}
\usepackage{graphicx} 
\usepackage[left=2cm, right=2cm, top=2cm, bottom=2cm]{geometry}
\usepackage{algorithm}
\usepackage{algpseudocode}
\usepackage{amsmath,amssymb}
\usepackage{comment}
\usepackage{hyperref}
\usepackage{stmaryrd}
\usepackage{bm}
\usepackage{authblk}
\usepackage{subcaption}
\usepackage{multirow}
\usepackage{booktabs} 

\graphicspath{{pic/}}
\title{Tensor-based compression of the sea temperature data}

\author[1,2]{Ilya Kosolapov}
\author[1,2]{Tatiana Sheloput}
\author[3,1]{Sergey Matveev\thanks{Corresponding author: matseralex@cs.msu.ru}}

\affil[1]{Marchuk Institute of Numerical Mathematics RAS, Moscow, Russia}
\affil[2]{Moscow Institute of Physics and Technology (National Research University), Dolgoprudnii, Russia}
\affil[3]{Lomonosov MSU, Faculty of Computational Mathematics and Cybernetics, Moscow, Russia}

\date{}

\begin{document}
\maketitle{}
\begin{abstract}
In this work we investigate efficient data compression for spatiotemporal Black, Azov and Marmara Seas temperature tensors that contain significant number of missing values. These tensors have a complex structure influenced by the coastlines and bathymetry, as well as temporal temperature changes. While such missing data typically provokes utilization of tensor completion algorithms, we demonstrate that standard SVD-based compression approaches (including the Tucker, Tensor-Train (TT) and Quantized-TT formats) are remarkably effective and yield comparable results. We propose a greedy spatial data partitioning algorithm enhancing their performance. We divide the data into the smaller subtensors before compression via exploitation of this trick.

Furthermore, our analysis reveals a strong temporal dependency in the data's compressibility caused by its nature. Fixing the level of precision we observe a significant seasonal variation. Investigating this, we find that a temporal partitioning on a scale of approximately two days is nearly optimal for all tested tensor based formats. The combined application of these spatial and temporal strategies with tensor methods ultimately achieves a robust compression ratio of 5 times across the entire dataset.            
\end{abstract}


\section{Introduction}

In modern Earth sciences, efficient data management is crucial for processing, storing, and transmitting the vast amounts of information generated by numerical models. This is particularly true in oceanography, where high-resolution spatiotemporal simulations produce data volumes that are challenging to handle. The temperature fields of the Black, Azov, and Marmara Seas, generated by the INM RAS model (see \cite{Zalesny,Fomin}) are studied in this paper as a prime example. These datasets form multidimensional tensors that describe complex physical processes, but their practical utility, as well as the ability to store and use them in tensor formats, is impacted by two main factors: their significant volume and the presence of extensive gaps, which arise from the irregular coastlines of the region.

Traditionally, the multidimensional data compression has been approached through two primary groups of methods: lossless compression preserves all information but offers limited compression ratios, and lossy compression that allows to control loss of information \cite{sullivan2004h, antsiferova2022video}. One can achieve much higher compression ratios via the lossy compression methods and it is often sufficient for analysis and visualization \cite{kurilovich2023case}. These methods are particularly suitable for multidimensional geophysical data fields, which typically exhibit high spatial and temporal correlations. This means that values at neighboring points in space or at subsequent time steps are usually close to each other, creating data redundancy analogous to that found in digital images \cite{Nelson, jia2022low}.

The modern data compression methods include tensor decompositions , such as the Tucker \cite{lieven2009survey, tucker1966some, sidiropoulos2017tensor}, Tensor-Train (TT) \cite{oseledets2009breaking, oseledets2011tensor}, and Quantized-TT (QTT) decompositions \cite{oseledets2009approximation}, effectively reveal hidden correlations and redundancy by representing the original tensor in a compact form (composed of a core tensor/tensors and factor matrices).

A common challenge in lossy compression of such datasets is handling missing data. The standard conceptual approach is to use tensor completion algorithms. To this end, we propose a compression framework based on the Singular Value Projection (SVP) algorithm \cite{meka2009guaranteed, petrov2023matrix, lebedeva2011tensor} in the Tucker format, that we use as a baseline for tensor reconstruction. This framework follows a scheme where only a small number of tensor entries, along with their mask $\Omega$, are stored. Reconstruction is then performed by applying the SVP algorithm to solve the completion problem and infer the missing entries, thereby recovering an approximation of the original tensor.

However, in this work we show a more direct and effective pathway. We first preprocess the incomplete dataset by dividing it into the set of regular, continuous subtensors using a greedy spatial partitioning algorithm. We compress these resulting subtensors compressed using the standard SVD-based tensor formats.  We show that this combination of partitioning and compression not only handles the original data gaps effectively but also yields superior results, outperforming the completion-based approach.

Furthermore, we analyze the data compressibility and demonstrate a strong temporal dependency, with significant seasonal variation. We show the comparison of data from January and May and obtain a degradation in performance for the Black, Azov, and Marmara Seas during the spring season.

The goal of this work is to systematically investigate our approach combining several data partitioning tricks with the established tensor decomposition frameworks. We evaluate the performance of Tucker, TT, and QTT decompositions on the sea temperature data and demonstrate the gains achieved through our spatial and temporal partitioning strategies, and assess the impact of the compression ratio on the accuracy of the reconstructed temperature fields.

\section{Dataset and its key features}
\label{sec:TT-method}
In order to analyze the efficiency and feasibility of the tensor-based compression algorithms we use the reanalysis dataset. The reanalysis data were obtained using variational assimilation of sea surface temperature (SST) \cite{Zakharova2021} observations from different satellites into the INMOM-based hydrothermodynamics model of the Black and Azov Seas for 2018 year \cite{agoshkov2018variational}. For assimilation, we exploit satellite SST data obtained from the Aqua, Terra, Suomi NPP, and Sentinel-3 spacecraft covering the Black, Azov, and Marmara Seas \cite{Agoshkov2009,Agoshkov2015,Parmuzin2012,Lupyan}.

The satellites provide data through distinct instruments with varying capabilities. The Terra and Aqua platforms are equipped with MODIS (Moderate Resolution Imaging Spectroradiometer) instruments, which have 36 spectral bands and are actively used to monitor sea surface temperature. The Suomi NPP satellite carries an advanced infrared radiometer designed for measuring surface temperatures at high latitudes. The Sentinel-3 satellite performs measurements with two primary instruments: the Sea and Land Surface Temperature Radiometer (SLSTR), which uses infrared cameras to accurately measure thermal radiation, and the Ocean and Land Color Instrument (OLCI), which can also provide SST data \cite{satellites_description}. It is important to note that satellite data constitutes instantaneous measurements that may not provide complete spatial coverage. Gaps can occur due to cloud cover, and some data may be filtered during processing due to quality issues. Moreover, the raw satellite input may contain outliers and isolated measurements where adjacent grid points have undefined values.

We construct a tensor $\mathcal{X}$ of temperature values basing on the specified reanalysis data. The position of a point on the sea surface $\Omega$ is described in a spherical coordinate system $(x, y, \sigma)$, or more precisely, by a pair of horizontal coordinates: $x$ (longitude) and $y$ (latitude). Further, we introduce regular discretizations in longitude ($x_1 < x_2 < \dots < x_N$) and latitude ($y_1 < y_2 < \dots < y_M$), such that for any point $\vec{p} = (x, y, \sigma(x, y))$ where $\vec{p} \in \Omega$, the conditions $x_1 \le x \le x_N$ and $y_1 \le y \le y_M$ are satisfied.

\begin{figure}[H]
\centering
\includegraphics[scale=0.6]{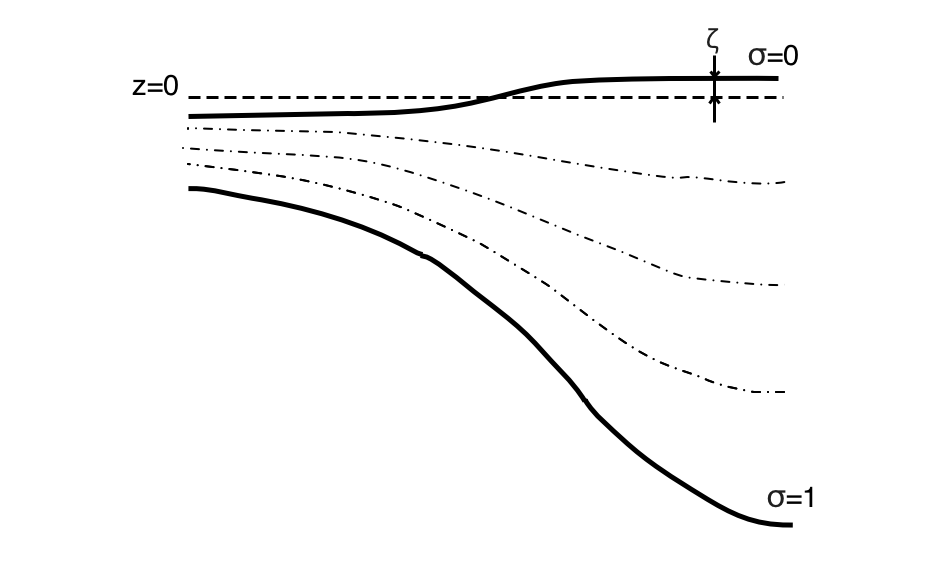}
\caption{Schematic distribution of $\sigma$-levels by ocean depth.}
\label{fig:sygma_system}
\end{figure}

The reanalysis data are available not only for the surface layer but also for depth in a $\sigma$-coordinate system \cite{Tolik}. The vertical variable is a dimensionless quantity defined as
$$\sigma(x,y,z,t) = \frac{z - \zeta(x,y,t)}{H(x,y) + \zeta(x,y,t)}\in [0,1],$$
where $x$ is longitude, $y$ is latitude, $z$ is depth, $H(x, y)$ is the ocean depth at rest, and $\zeta(x, y, t)$ is the sea level deviation from the unperturbed surface (see Fig. \ref{fig:sygma_system}). In this work, we consider the a system of $L$ fixed vertical levels $\sigma_1 < \ldots < \sigma_L$, where the sea level deviation is assumed to be identically zero: $\zeta(x, y, t) \equiv 0$.

In addition to spatial coordinates, the sea temperature also depends on the time $t$. Thus, a complete description of the temperature field requires the introduction of a time grid $t_1 < t_2 < \dots < t_K$.

As a result, the reanalysis data can be represented as a four-dimensional regular array (tensor) $\mathcal X$ of size $N \times M \times L \times K$, where:
\begin{itemize}
    \item $N$ is the number of grid nodes in longitude,
    \item $M$ is the number of grid nodes in latitude,
    \item $L$ is the number of depth levels,
    \item $K$ is the number of time slices.
\end{itemize}

An element of this tensor $\mathcal{X}_{i,j,l,k}$ represents the temperature value at the point with coordinates $(x_i, y_j)$ at depth $\sigma_l$ at time $t_k$. Unlike the original satellite data, this array is regular and complete over the ocean domain. The temperature value in our data is defined for every valid combination of indices (i.e., where the coordinate $(x_i, y_j)$ corresponds to water). If the coordinate $(x_i, y_j)$ does not belong to water, the corresponding values of $\mathcal{X}_{i,j,l,k}$ are marked as missing values.

In our numerical experiments, we store the temperature data on a rectangular grid of $306 \times 200$ points with the following physical dimensions:
\begin{itemize}
\item \textbf{Longitude:} 306 samples with a step of 0.05°
\item \textbf{Latitude:} 200 samples with a step of 0.036°
\item \textbf{Depth:} 20 $\sigma$-levels
\item \textbf{Time:} 2922 time slices for the year 2018 with a uniform step of 3 hours
\end{itemize}

Summing up, in our numerical experiments we used the tensor $\mathcal X \in \mathbb R^{N \times M \times L \times K}$ of dimension $(306, 200, 20, 2922)$. In some tests, we studied the temporal subsamples resulting in the tensor dimension $(306, 200, 20, T)$,
where $T$ is the length of a sequential time segment.

\section{Notations and methods}

\subsection{Tucker decomposition}

In this section we discuss the necessary information about the tensors and Tucker decomposition (see \cite{tucker1966some, kolda2009tensor, lieven2009survey}).
A tensor of dimension $d$ or a $d$-dimensional tensor is $\mathcal{X} \in \mathbb{R}^{n_1 \times \cdots \times n_d}$. An element of the tensor $\mathcal{X}$ at position $(i_1, \dots, i_d)$ is denoted as $\mathcal{X}(i_1, \dots, i_d)$ or $x_{i_1 \dots i_d}$. The columns of the tensor in the $k$-th direction are vectors
\[
\mathcal{X}(i_1, \dots, i_{k-1}, :, i_{k+1}, \dots, i_d) \in \mathbb{R}^{n_k},
\]
which are obtained by fixing all indices except the $k$-th.

\textit{The unfolding} in the $k$-th direction of the tensor $\mathcal{X} \in \mathbb{R}^{n_1 \times \cdots \times n_d}$ is the matrix $\mathcal{X}_{(k)} \in \mathbb{R}^{n_k \times \prod_{i \neq k} n_i}$ that is composed of columns in the $k$-th direction, ordered in such a way that firstly the index $i_1$ changes, then the index $i_2$ and so on. For tensors of the fixed size we can introduce a scalar product and Frobenius norm generated by this scalar product
\[
\langle \mathcal{X}, \mathcal{Y} \rangle_F = \sum_{i_1, \dots, i_d} (\mathcal{X} \circ \mathcal{Y})_{i_1 \dots i_d}, \quad \|\mathcal{X}\|_F = \sqrt{\langle \mathcal{X}, \mathcal{X} \rangle_F},
\]
where $\circ$~--- is the Hadamard (element-wise) product. This scalar product does not depend on the order of the elements, so
\[
\langle \mathcal{X}, \mathcal{Y} \rangle_F = \langle \mathcal{X}_{(k)}, \mathcal{Y}_{(k)} \rangle_F.
\]
In our work we also have to deal with the Chebyshev norm
\[\|\mathcal X\|_C = \max_{i_1,\ldots, i_d} x_{i_1,\ldots, i_d}\]
because we need to measure the largest deviation in the approximation from the original sea temperature data in the usual Celsius degrees. We define the multiplication operation of the tensor $\mathcal{X} \in \mathbb{R}^{n_1 \times \cdots \times n_d}$ by the matrix $A \in \mathbb{R}^{m_k \times n_k}$ in the $k$-th direction, which is denoted by $\mathcal{X} \times_k A$. Its result is the tensor $\mathcal{Y}$:
\[
\mathcal{Y} = \mathcal{X} \times_k A \in \mathbb{R}^{n_1 \times \cdots \times n_{k-1} \times m_k \times n_{k+1} \times \cdots \times n_d}, \quad \mathcal{Y}_{(k)} = A \mathcal{X}_{(k)}.
\]
The products of the tensor and matrices in different directions commute:
\[
\mathcal{X} \times_k A \times_l B = \mathcal{X} \times_l B \times_k A, \quad \text{if} \quad k \neq l.
\]
One may do it also along a single direction:
\[
\mathcal{X} \times_k A \times_k B = \mathcal{X} \times_k (BA).
\]

In these notations one can define the \textit{Tucker decomposition} of a tensor $\mathcal{X} \in \mathbb{R}^{n_1 \times \cdots \times n_d}$ as following
\[
\mathcal{X} = \mathcal{G} \times_1 U_1 \times_2 \cdots \times_d U_d, \quad \mathcal{G} \in \mathbb{R}^{r_1 \times \cdots \times r_d}, \quad U_k \in \mathbb{R}^{n_k \times r_k},
\]
where the tensor $\mathcal{G}$~-- \textit{decomposition kernel}, matrices $U_k$~-- \textit{factors}, $r_k$~-- \textit{decomposition ranks}.

If for some elementwise precise decomposition of the tensor $\mathcal{X}$ one may prove that all its ranks $r_k$ are minimal among all its Tucker decompositions, then it is called \textit{minimal} decomposition, and the $r_k$ are called \textit{Tucker ranks} (or multilinear ranks) \cite{lieven2009survey}. The minimal Tucker decomposition always exists and the Tucker ranks are equal to the usual matrix ranks of the unfoldings in the corresponding directions
\[
r_k = \operatorname{rank}(\mathcal{X}_{(k)}).
\]
\textit{Tucker rank} of a tensor $\mathcal{X}$ is an ordered set of its ranks
\[
\operatorname{trank}(\mathcal{X}) = (r_1, \dots, r_d).
\]
If $\operatorname{rank}(U_k) = r_k$ for all $k$, then the Tucker ranks of the tensor $\mathcal{X}$ and the kernel $\mathcal{G}$ coincide: $\operatorname{trank}(\mathcal{X}) = \operatorname{trank}(\mathcal{G})$.

If the factors of the Tucker decomposition $U_k$ have orthogonal columns, then the decomposition is called \textit{orthogonal}. For an orthogonal decomposition, it is true
\[
\|\mathcal{X}\|_F = \|\mathcal{G}\|_F.
\]

The famous Higher-Order Singular Value Decomposition (HOSVD) algorithm constructs a minimal (orthogonal) Tucker decomposition \cite{lieven2009survey}. This method can be considered as a tensor-based analogy of the singular value decomposition (SVD) for matrices \cite{golub2013matrix, tyrtyshnikov2012brief}.

\begin{algorithm}
\begin{algorithmic}[1]
\Function{HOSVD}{$\mathcal{X}, r_1, r_2, \ldots, r_N$} 
\For{$k = 1$ \textbf{to} $N$} 
\State $A^{(k)} \gets$ $r_k$ leading left singular vectors $\mathcal{X}_{(k)}$ 
\EndFor 
\State $\mathcal{G} \gets \mathcal{X} \times_1 A^{(1)^T} \times_2 A^{(2)^T} \cdots \times_N A^{(N)^T}$ 
\State \Return $\mathcal{G}, A^{(1)}, A^{(2)}, \ldots, A^{(N)}$
\EndFunction
\end{algorithmic}
\caption{Pseudocode of the Higher-Order Singular Value Decomposition.}
\label{algo:HOSVD}
\end{algorithm}
In contrast to the SVD for matrices, the HOSVD in application to tensors is not exactly an \textit{optimal projection} of $\mathcal{P}_r$ onto the set with the fixed ranks $\mathcal M_r$
\[\| \mathcal{P}_r(\mathcal{X}) - \mathcal{X} \|_F = \min_{\hat{\mathcal{X}}\in Mr}\|\hat{\mathcal{X}} - \mathcal{X}\|_F,\]
but fortunately has the quasi-optimality property
\[
\|\operatorname{HOSVD}_r(\mathcal{X}) - \mathcal{X}\|_F \leq \sqrt{d} \|\mathcal{P}_r(\mathcal{X}) - \mathcal{X}\|_F,
\]
where $\mathcal{P}_r$ is an optimal projection onto the set tensors of Tucker rank at most $r$. This prorerty is extensively used for theoretical justification of various numerical algorithms storing the data in Tucker format \cite{sultonov2023low, TuckER, Tao}.

\subsection{Tensor completion in Tucker format}

Quite often one may meet the situation where some observations are missing during processing of the multidimensional data  \cite{TuckER, Tao}. These elements of the data may be unavailable due to the hardware failures \cite{petrov2023theoretical, petrov2024constructing}, measurement limitations \cite{sheloput2024application}, or many other factors \cite{petrov2023matrix, budzinskiy2023tensor, frolov2017tensor, lebedeva2011tensor}.  The general task of the matrix and tensor completion methods is to reconstruct these missing elements basing on the available observations. Assume that for some tensor $\mathcal{X} \in \mathbb{R}^{n_1 \times n_2 \times \cdots \times n_d}$ we have
\begin{enumerate}
\item The set of positions of known elements: $\Omega \subset \{(i_1, i_2, \dots, i_d) \mid 1 \leq i_k \leq n_k\}$
\item Projection operator $\mathcal{P}_\Omega$ defining the available elements
\[
\mathcal{P}_\Omega(\mathcal{X})_{i_1\dots i_d} =
\begin{cases}
t_{i_1\dots i_d}, & (i_1,\dots,i_d) \in \Omega \\
0, & \text{otherwise.}
\end{cases}
\]
\end{enumerate}
In our case, the sea temperature array contains for a lot of formally missing points. These elements arise due to natural inhomogeneity of the seabed as well as its irregular spatial form that is not a rectangle. Thus, we need to find the full tensor $\mathcal{\hat X} \in \mathbb{R}^{n_1 \times \cdots \times n_d}$, satisfying the two properties:
\begin{enumerate}
\item it matches with known data
\[
\mathcal{P}_\Omega(\mathcal{\hat X}) = \mathcal{P}_\Omega(\mathcal{X}),
\]
\item has the low-rank stucture allowing to compress the data
\[
\text{trank}(\mathcal{\hat X}) = (r_1, \dots, r_d) \quad \text{c} \quad r_k \ll n_k \quad \forall k,
\]
where $\operatorname{trank}$ are the Tucker ranks.
\end{enumerate}

The problem of completion in the Tucker format can be formulated as a rank optimization problem with constraints
\[
\min_{\mathcal{\hat X}} \ \operatorname{trank}(\mathcal{\hat X}) \quad \text{subject to} \quad \mathcal{P}_\Omega(\mathcal{\hat X}) = \mathcal{P}_\Omega(\mathcal{X}),
\]
where $\operatorname{trank}(\mathcal{\hat X})$ is the Tucker rank of the tensor $\mathcal{\hat X}$.

Since the direct finding the minimal rank and decomposition is a hard problem, we replace it by a simpler one
\[
\min_{\mathcal{\hat X}} \ \|\mathcal{P}_\Omega(\mathcal{\hat X})- \mathcal P_\Omega(\mathcal X)\|_F \quad \text{subject to} \quad \operatorname{trank}(\mathcal{\hat X})\leq r.
\]
This approach allows one to search for the tensor in factorized form and use gradient methods. This problem can be solved via various approaches \cite{petrov2023matrix, budzinskiy2023tensor, frolov2017tensor}. For example, one may use the SVP-typed \cite{meka2009guaranteed, petrov2024constructing} (singular value projection) algorithm \ref{algo:SVP}, where the rank update procedure is $r_k = \min(\hat r_k, r+\delta)$, $\delta$ is the rank increase constant at each step.
\begin{algorithm}
\begin{algorithmic}[1]
\Function{SVP}{$\mathcal{P}_{\Omega}(\mathcal X),\hat{r} ,\varepsilon , \delta, \eta_t$}
\Comment{$\mathcal{P}_{\Omega}(\mathcal X)$ - observed elements, $\hat{r}$ - target Tucker rank, $\varepsilon$ - relative error, $\delta$ -- rank update parameter,
$\eta_t$ - gradient step}
\State $\hat{\mathcal X} = \mathcal{G} \times_1 U_1 \ldots \times_d U_d \leftarrow \mathbf{0}$ \Comment{Zero tensor}
\State $r \leftarrow [1, 1, \dots, 1]$ \Comment{Initial ranks} 

\While{$\dfrac{\|\mathcal{P}_{\Omega}(\mathcal{\hat X} - \mathcal{X})\|_F}{\|\mathcal{P}_{\Omega}(\mathcal{X})\|_F} \geq \varepsilon$} 
\State $\mathcal{\hat X} \leftarrow \mathcal{\hat X} - \eta_t \mathcal{P}_{\Omega}(\mathcal{\hat X} - \mathcal{X})$ 
\Comment{Gradient step for elements on the mask} 
\State $\mathcal{G} \times_1 U_1 \ldots \times_d U_d \ \leftarrow \text{HOSVD}_r(\mathcal{\hat{X}})$ \Comment{Truncation to rank $r$ of the whole tensor} 
\State $\hat{r} \leftarrow \text{update\_rank}(r, \delta)$ \Comment{Update Tuckers ranks with shift parameter $\delta$}
\EndWhile
\State \Return $\mathcal{\hat{X}}$
\EndFunction 
\end{algorithmic}
\caption{Pseudocode of the SVP-typed procedure in the Tucker format.}
\label{algo:SVP}
\end{algorithm}

\subsection{Tensor-Train decomposition}

Let us start with definition of the tensor train for the three-dimensional data. In this case, the tensor train corresponds to the following representation
\begin{align*}
A(i,j,k) = \sum_{\alpha_1}^{r_1} \sum_{\alpha_2}^{r_2} G_1(i, \alpha_1) G_2(\alpha_1, j, \alpha_2) G_3(\alpha_2, k).
\end{align*}
In some sense, this representation can be considered a ``not yet finalized'' Tucker decomposition, requiring $N r_1 r_2 + N r_1 + N r_2$ memory cells. At first glance, not seeming like a good idea for data compression. However, the remarkable properties of this format become apparent in the case of high-dimensional data, since it does not explicitly contain factors that have exponential memory requirements with increasing dimension $d$ as we may find for the Tucker decomposition.

Fortunately, this decomposition can also be calculated using a procedure based on the singular value decomposition. This procedure is known as Tensor Train Singular Value Decomposition (TTSVD) \cite{oseledets2011tensor} and is written as follows for three-dimensional arrays:
\begin{align*}
& A(i,\overline{i,k}) = \sum_{\alpha_1=1}^{r_1} G_1(i, \alpha_1) G_{2,3}(\alpha_1, \overline{j,k}) \quad \{\text{use the SVD with rank } r_1 \text{ for } A(i,\overline{j,k})\}
\\
& = \sum_{\alpha_1=1}^{r_1} G_1(i, \alpha_1) G_{2,3}(\overline{\alpha_1, j}, k) \quad \{\text{use the reshape for } G_{2,3}(\alpha_1, \overline{j,k}) \text{ and SVD with rank } r_2 \text{ for } G_{2,3}(\overline{\alpha_1, j},k)\}\\
& = \sum_{\alpha_1=1}^{r_1} \sum_{\alpha_2=1}^{r_2} G_1(i, \alpha_1) G_{2}(\alpha_1, j, \alpha_2) G_3(\alpha_2, k).
\end{align*}
From the steps above, it is clear that the complexity of obtaining such a decomposition is $O(N^4)$ operations. The minimal possible values of the TT ranks correspond to the ranks of the unfolding matrices (slightly different from those in the Tucker decomposition):
\begin{equation*}
A_1 = A(i, \overline{j,k}) \text{ and } A_2 = A(\overline{i,j}, k).
\end{equation*}
The full the tensor train in the general d-dimensional case corresponds to the following representation of original array
\begin{equation}\label{eq:tensor_train}
A(i_1,i_2,\ldots, i_d) = \sum_{\alpha_1}^{r_1} \sum_{\alpha_2}^{r_2} \ldots \sum_{\alpha_{d-1}}^{r_{d-1}} G_1(i_1, \alpha_1) G_2(\alpha_1, i_2, \alpha_2) \ldots G_{d}(\alpha_{d-1}, i_d)
\end{equation}
and can be constructed in the same manner via the general TTSVD procedure \cite{oseledets2011tensor} within $O(N^{d+1})$ operations. The final representation requires just $O(d N r^2)$ memory cells instead of original $N^d$, where $r=\underset{i=1~, \ldots,~n}{\max} r_i$. It is easy to see that such a representation exists for any $d$-dimensional array (maybe with high ranks) and in other literature souces it is known as a matrix product state. This name seems quite reasonable, since computing an element of the array $A(i_1,i_2,\ldots, i_d)$ represented in the tensor train format can be accomplished via $(d-2)$ matrix-vector multiplication and one final dot product.

If $r=\underset{i=1~, \ldots,~n}{\max} r_i$ is small then such a representation leads to the fantastic compression of data. The minimal ranks of the d-dimensional tensor train are equal to the ranks of the unfolding matrices
\begin{equation*}
A_1 = A(i_1, \overline{i_2, \ldots, i_k}),~ A_2 = A(\overline{i_1,i_2}, \overline{i_3\, \ldots, i_d}),~ \ldots, A_{d-1} = A(\overline{i_1,\ldots,i_{d-1}}, i_d),
\end{equation*}
and can be formally found via the SVD. Fortunately, the approximate rank truncation via TTSVD has the similar quasi-optimality properties to HOSVD. Hence, in the approximations we obtain that
$$||A - B_t||_F \leq \sqrt{d-1}||A - B_{*}||_F,$$
where $B_t$ is the result of TTSVD, and $B_{*}$ is the best possible approximation in TT format with the same fixed ranks. 

In addition to the broad family of low-rank tensors with exact rank values, there is a broader family of tensors \cite{tyrtyshnikov2003tensor} with very good approximations via the low-rank tensor train (e.g., in the Frobenius norm) that can be computed either using the quasioptimal TTSVD algorithm (and its randomized generalizations) \cite{oseledets2009breaking, sultonov2023low}, or using the much faster but less accurate adaptive TT-cross approximation approach \cite{oseledets2010tt, oseledets2011tensor} requiring the computation of just $O(d N r^3)$ elements of the target tensor. Although this format has been known for a long time as the matrix product state \cite{white1992density, verstraete2023density} but its independent reinterpretation using numerical linear algebra methods gave rise to algorithms leading to constructive approximations of the miscellaneous data.

Currently, there are many useful and efficient basic methods (such as element-wise sum and product) \cite{khoromskij2018tensor} for manipulating with data in the tensor train format, as well as sophisticated optimization procedures such as AMEn (a method based on minimizing the energy functional) for high-dimensional linear systems \cite{dolgov2014alternating} or heuristic global optimization procedures \cite{zheltkov2020global}.

The tensor train representations has also shown drastically high performance for data compression via utilization of the virtual dimensions even into usual large vectors. This idea allows one to transform an array with $N=2^d$ elements into a d-dimensional tensor $\underbrace{2 \times \ldots \times 2}_{\text{d times}}$ that ca be further decomposed as a tensor train. If the ranks of the corresponding decomposition (either exact or approximate) then one may see a colossal compression of the original data to $O(d R^2)$ memory cells meaning the logarithmic requirements of $O(R^2 \cdot \log N)$ elements for storing the data. This trick is known as the quantized tensor train (QTT) format \cite{oseledets2009approximation}. Moreover, one may also introduce the virtual dimensions reshaping the original modes of the multidimensional arrays (see e.g. \cite{danis2025tensor, novikov2015tensorizing}) increasing the level of their compression within the QTT format. We also try this idea in application to our data below in our work.

\section{Greedy partitioning of data into sub-tensors}

The output of the reanalysis model is the structured four-dimensional tensor $\mathcal{X}$ with missing values. We obtain a significant challenge with its structured representation because this regular grid encompasses both ocean and land masses. The temperature values are naturally undefined at the grid nodes corresponding to the land creating a sparse data structure within a dense array. Applying standard tensor decompositions in  to $\mathcal{X}$ straight-forward manner is inefficient and even mathematically unsound as these operations would involve invalid land data.

In order to enable the efficient computing we introduce a strategy that decomposes the full tensor into a collection of smaller, continuous regular sub-tensors that we refer to as \textit{blocks}. We denote each block as $\mathcal X_I\in \mathbb R^{(x_{\text{end}} - x_{\text{start}}) \times (y_{\text{end}}-y_{\text{start}})\times K\times T}$, is defined by its horizontal indices $I=(x_{\text{start}}, x_{\text{end}}, y_{\text{start}},y_{\text{end}})$. Hence, it contains only the valid ocean data. The vertical ($L$) and temporal ($T$) dimensions are included in full for any valid horizontal region.

A critical observation is that the pattern of defined (ocean) and undefined (land) data is consistent across all depth levels and time steps. It is determined solely by the coastline. Let us define the \textit{$j$-th layer of a tensor} $\mathcal X^{j}$ as the section of the tensor in the direction $\sigma$, i.e. $\mathcal X^j=\mathcal X[:,:,j,:]$. 

Therefore, the problem of partitioning the 4D tensor reduces to finding a set of valid rectangular blocks on a single 2D horizontal slice (e.g., the surface layer $\mathcal{X}^1$ at any time). The resulting block structure should be further applied uniformly across all depth levels and time steps. The block is considered valid if it does not contain any land points. Hence, all values in the corresponding region of the data mask should be defined.

\begin{figure}[H]
\centering
\includegraphics[width=0.8\linewidth]{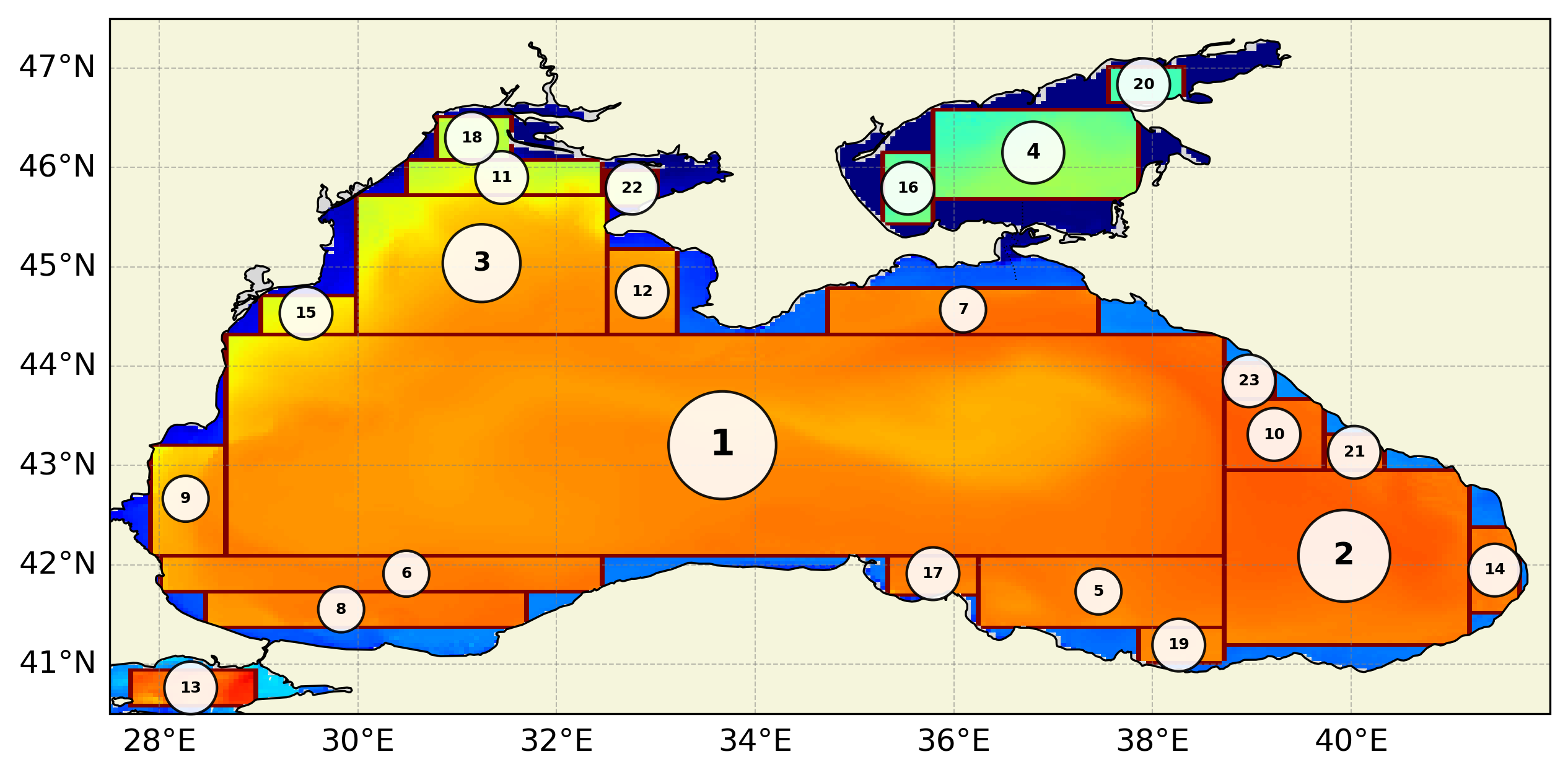}
\caption{Results of greedy partitioning of sea temperature data layer into blocks. In total we have 23 blocks}
\label{fig:greedy_hmap}
\end{figure}

All, in all, we propose a greedy algorithm to partitioning of the 2D surface data into the valid rectangular blocks of the maximum possible size. The algorithm iteratively finds the largest possible valid block in the unallocated area of the domain, marks that area as used, and repeats until the entire domain is covered. The algorithm requires a minimum block size $S_{\min}$ to prevent the creation of excessively small blocks. The pseudocode for the main procedure and its subroutines is provided in Algorithms~\ref{alg:greedy}, \ref{alg:largest_block}, and \ref{alg:IsValidBlock}. Figure~\ref{fig:greedy_hmap} demonstrates the application of this algorithm to our dataset, while Table~\ref{tab:block_dimensions} summarizes the sizes of the resulting subtensors derived from spatial partitioning. 

Our idea of data partitioning seems to be coherent with idea of the wavelet tensor train decomposition \cite{kharyuk2014wtt} that has been also applied to compression of the multidimensional fMRI data \cite{pasha2014compression}. In contrast to these works our greedy approach does not utilize the wavelet transformation for for better structuring of the studied data.

\begin{table}[H]
    \centering
     {\footnotesize
    \begin{tabular}{|c|c|c|c|c|c|c|c|c|c|c|c|c|}
        \hline
        № & 1 & 2 & 3 & 4 & 5 & 6 & 7 & 8 & 9 & 10 & 11 & 12 \\
        \hline
        Modes & $199 \times 62$ & $49 \times 49$ & $50 \times 39$ & $41 \times 25$ & $49 \times 20$ & $88 \times 10$ & $54 \times 13$ & $64 \times 10$ & $15 \times 31$ & $20 \times 20$ & $39 \times 10$ & $14 \times 24$ \\
        \hline
    \end{tabular}

\vspace{0.3cm}

    \begin{tabular}{|c|c|c|c|c|c|c|c|c|c|c|c|}
        \hline
        № & 13 & 14 & 15 & 16 & 17 & 18 & 19 & 20 & 21 & 22 & 23 \\
        \hline
        Modes & $25 \times 10$ & $10 \times 24$ & $19 \times 11$ & $10 \times 20$ & $18 \times 11$ & $15 \times 12$ & $17 \times 10$ & $15 \times 10$ & $12 \times 10$ & $10 \times 10$ & $10 \times 10$ \\
        \hline
    \end{tabular}
    
    \caption{Horizontal dimensions of the partitioned blocks in mode format. The first mode represents latitude and the second mode represents longitude. Each block contains 20 vertical layers in the $\sigma$-coordinate system. The fourth dimension (time) represents the number of subsequent time steps $T$. The overall dimensions of each block are structured as latitude $\times$ longitude $\times$ 20 $\times$ $T$.}
    \label{tab:block_dimensions}
    }
\end{table}

\begin{algorithm}[H]
\begin{algorithmic}[1]
\Function{GreedyPartitioning}{$\mathcal X^1, S_{\min}$}
\Comment{$S_{\min}$ - minimum block size}
\State $M \leftarrow \mathbf{false}$
\Comment{Mask of already used blocks}
\While{True}
\State $I\leftarrow \text{FindLargestBlock}(M)$
\If{$I = \text{None}$}
\State \text{break}
\EndIf
\State $M_I \leftarrow 1$
\Comment{Filling mask by indices $I$}
\EndWhile
\State \Return $[I_1,\ldots, I_d]$
\EndFunction 
\end{algorithmic}
\caption{Pseudocode of the greedy data partitioning into the regular subtensors. }
\label{alg:greedy}
\end{algorithm}

\begin{algorithm}[H]
\begin{algorithmic}[1]
\Function{FindLargestBlock}{$M$}
\Comment{Mask of already used blocks}
\State $A_{\text{max}}\leftarrow0$
\Comment{Maximum block size}
\State $I_{\text{best}} = \text{None}$
\For{{$i = 0$ \textbf{to} $N_x - S_{\text{min}}$}}
\For{{$j = 0$ \textbf{to} $N_y - S_{\text{min}}$}}
\State $(x_{\min}, x_{\max}, y_{\min}, y_{\max}) \leftarrow (i, i + S_{\min},j, j + S_{\min})$ 
\If {\textbf{not} $\text{ IsValidBlock}(x_{\min}, x_{\max}, y_{\min}, y_{\max})$} 
\State\text{continue} 
\EndIf 
\While{$x_{\max} < N_x \textbf{ and } \text{IsValidBlock}(M,x_{\min}, x_{\max}+1, y_{\min}, y_{\max})$} 
\State $x_{\max} \leftarrow x_{\max} +1$ 
\EndWhile 
\While{$y_{\max} < N_y \textbf{ and } \text{IsValidBlock}(M,x_{\min}, x_{\max}, y_{\min}, y_{\max}+1)$} 
\State $y_{\max} \leftarrow y_{\max} +1$ 
\EndWhile 
\State $A \leftarrow (x_{\max}-x_{\min})\cdot(y_{\max}-y_{\min})$ 
\If{$A> A_{\max}$} 
\State $A_{\max} \leftarrow A$ 
\State $I_{\text{best}} \leftarrow (x_{\min},x_{\max},y_{\min},y_{\max})$ 
\EndIf 
\EndFor 
\EndFor 
\State \Return $I_{\text{best}}$
\EndFunction
\end{algorithmic}
\caption{Pseudocode finding the largest block during the greedy partitioning step.}
\label{alg:largest_block}
\end{algorithm}

\begin{algorithm}[H]
\begin{algorithmic}[1]
\Function{IsValidBlock}{$M,x_{\min},x_{\max},y_{\min},y_{\max}$}
\If{$x_{\max}-x_{\min} < S_{\min} \textbf{ or } y_{\max}-y_{\min} < S_{\min}$}
\State \Return \text{false}
\EndIf
\State $I \leftarrow (x_{\min},x_{\max},y_{\min},y_{\max})$
\State \Return $\text{any}(\mathcal X_I \neq \text{None}) \textbf{ and } \text{any}(\mathcal M_I \neq 1)$
\EndFunction 
\end{algorithmic}
\caption{Function that checks the block for being correct}
\label{alg:IsValidBlock}
\end{algorithm}

\section{Numerical experiments}

In this section we discuss the results of numerical experiments with application of the SVP-based completion in the Tucker format and also SVD-based thresholding in the Tucker and Tensor Train representations for compression of our data. In all experiments we measure the relative error of approximation in Frobenius norm as 
$$\frac{||\mathcal{X} - \hat{\mathcal{X}}||_F}{||\mathcal{X}||_F}.$$
Utilization of this criteria allows us to reach really great compression rates (by dozens or even hundreds of times as we show in Tables ) of the largest block in our partitioned data setting the level of accuracy of our approximation as $\varepsilon = 10^{-3}$. At the same time the absolute error of approximation 
$$\|\mathcal{X - \hat{\mathcal{X}}}\|_C$$
does not seem to be acceptable for us because it exceeds 2 or even 3 degree Celsius. Such large data approximation errors make them unreliable for further use in geophysical applications. The natural solution of this problem is increase of the ranks of the decomposition leading to much more modest levels of compression of our data (just several times). However, even its compression by 1.9-5 times with acceptable level of accuracy (less than 0.5 degree Celsius) is useful for economic storage of our data.

To evaluate the performance of tensor compression algorithms, we define two complementary compression ratios that provide transparent comparison metrics. Let $N(\mathcal{X})$ denote the number of elements required to store tensor $\mathcal{X}$, and let $\hat{\mathcal{X}}_I$ represent the compressed version of subtensor $\mathcal{X}_I$:

\[
CR_{\text{all}} = \frac{N(\mathcal{X})}{\sum_I N(\hat{\mathcal{X}}_I)}, \quad CR_{\text{sub}} = \frac{\sum_I N(\mathcal{X}_I)}{\sum_I N(\hat{\mathcal{X}}_I)}.
\]

The \textit{Overall Compression Ratio} ($CR_{\text{all}}$) measures the total compression achieved across the entire dataset. The \textit{Subtensor Compression Ratio} ($CR_{\text{sub}}$) isolates the compression efficiency within the spatial partitioning scheme itself, accounting for potential non-optimalities -- such as uncompressed data near coastlines. In other words, the uncompressed values are not taken into account when calculating $CR_{\text{sub}}$.

\subsection{SVP-based completion in the Tucker format}

At first, we designed the numerical experiment allowing to evaluate the efficiency of the Singular Value Projection (SVP) inspired algorithm \ref{algo:SVP} for data compression via tensor completion.  The matrix completion methodology has already been applied to compression of the sea surface temperature \cite{sheloput2024application} and one may infer that its generalization might work well for a case of the higher dimensionality.

\begin{figure}[H]
    \centering
    \includegraphics[width=0.6\linewidth]{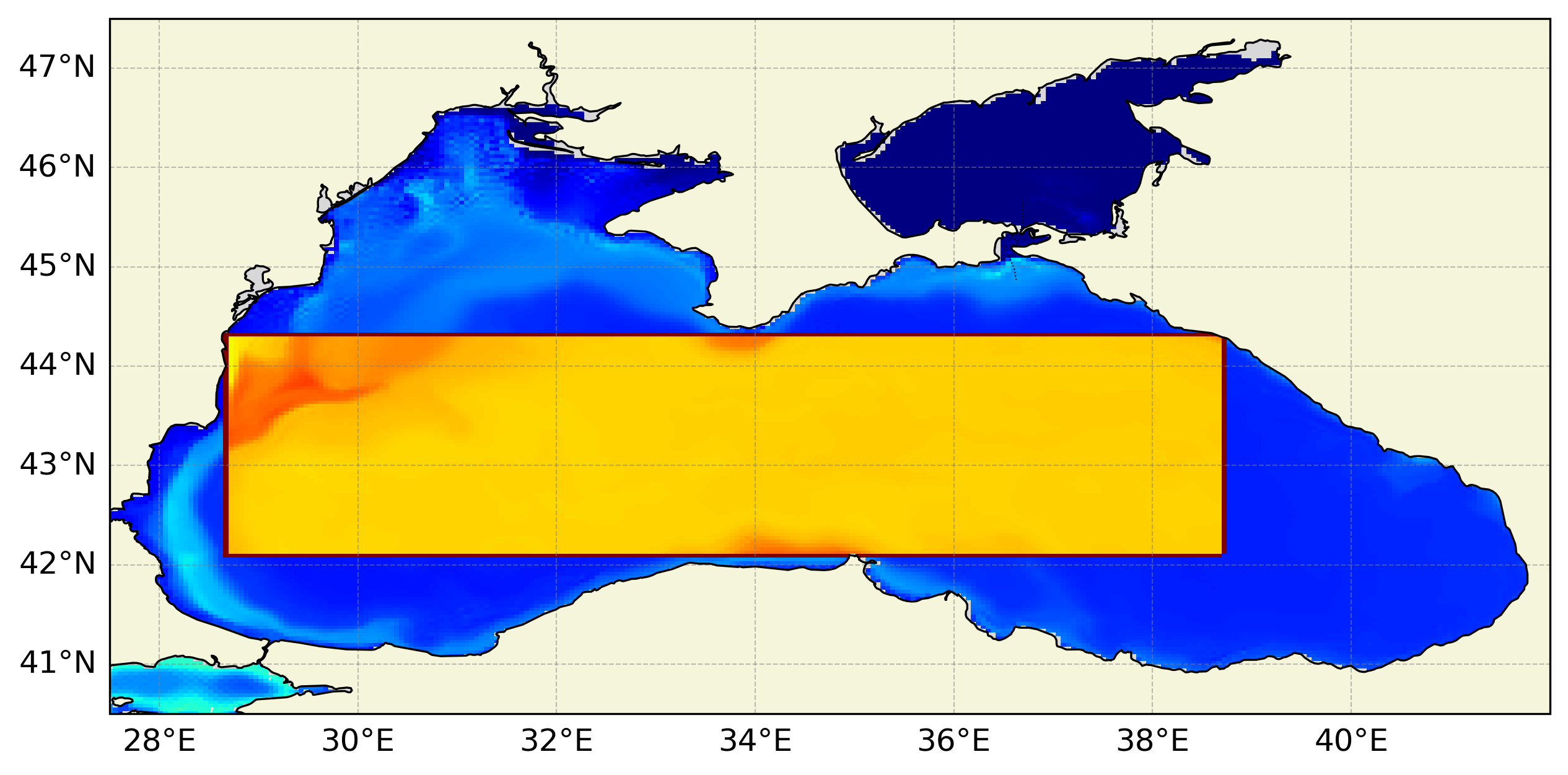}
    \caption{The largest block $\mathcal{X}_{I{\text{max}}}$}
    \label{fig:largest_block}
\end{figure}

For simplicity and clarity of analysis, we isolate the largest continuous block within the original dataset, denoted as $\mathcal{X}_{I{\text{max}}}$ (Figure \ref{fig:largest_block}). Further we apply a series of random sampling masks $\Omega$ to this tensor simulating the compression. We generate each mask $\Omega$ by random selection of a subset of the tensor entries to be ``observed''. We set the number of selected entries allowing to achieve the specific Compression Ratios (CR). 

Therefore, a sparser mask $\Omega$ corresponds to a higher CR, as fewer elements are retained. The core idea of the compression scheme is to store only these randomly selected entries and the mask $\Omega$ itself. After it we check the reconstruction results by using the SVP algorithm. We solve the completion problem numerically and infer that the missing entries have to be recovered by an approximation $\mathcal{\hat{X}}_{I{\text{max}}}$.

\begin{figure}[H]
    \centering
    \includegraphics[width=0.95\linewidth]{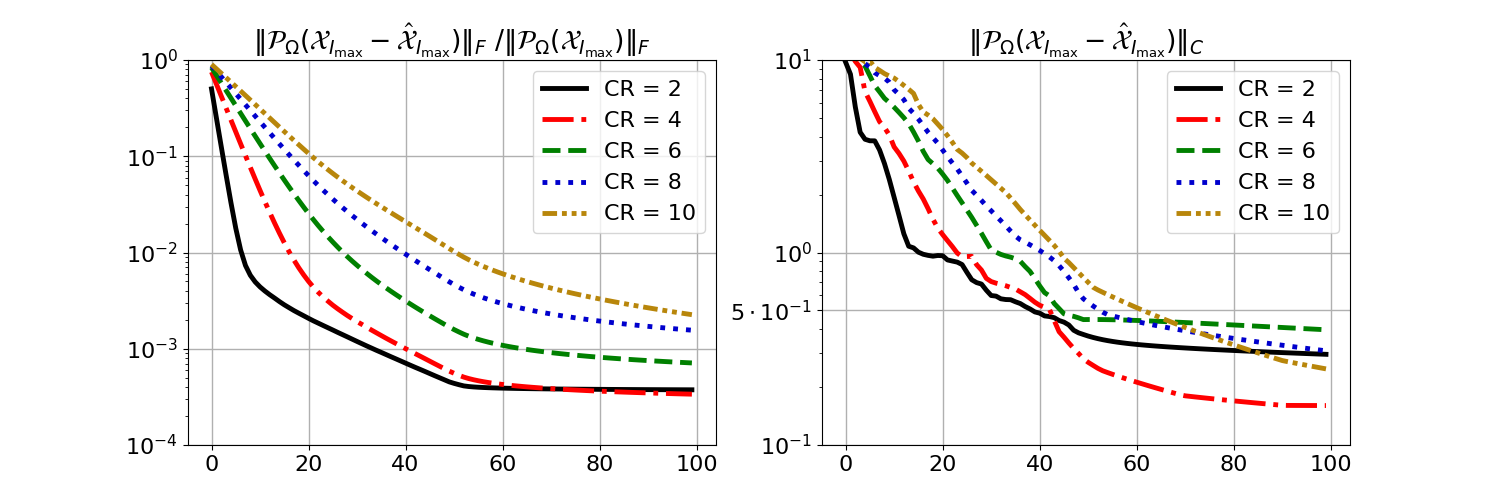}
    \caption{Convergence of the SVP algorithm on the random mask of the measured elements, assessed by the Frobenius norm (left) and Chebyshev norm (right) for $\mathcal{X}_{I_{\max}}$ of January data. The dependence on the number of iterations is shown for different compression ratios $CR = 2, 4, \dots, 10$.}
    \label{fig:SVP_convergence_jan}
\end{figure}

In figure \ref{fig:SVP_convergence_jan} we show the convergence characteristics of the Singular Value Projection (SVP) algorithm for the low-rank tensor completion task, specifically applied to the January dataset $\mathcal{X}_{I{\text{max}}}$. The convergence is evaluated using two  metrics: the Frobenius norm (left panel) that measures the overall global error, and the Chebyshev (infinity) norm (right panel), which captures the magnitude of the largest individual error element.

Each trajectory in the plot corresponds to a different value of compression ratio (CR) from 2 to 10. The Frobenius norm error decreases monotonically with each iteration as we may guess from the structure of the completion procedure. 

The analysis of the Chebyshev norm reveals a more nuanced behavior. While the initial errors are larger for higher CRs, their trajectories eventually cross those of lower CRs, stabilizing at a lower final error level. Hence, we may assume that the random sampling mask $\mathcal{P}_{\Omega}(\cdot)$, while introducing high initial uncertainty, may ultimately allow the algorithm to find a more robust solution that avoids large outliers across the approximated tensor.

\begin{figure}[h]
    \centering
    \includegraphics[width=0.95\linewidth]{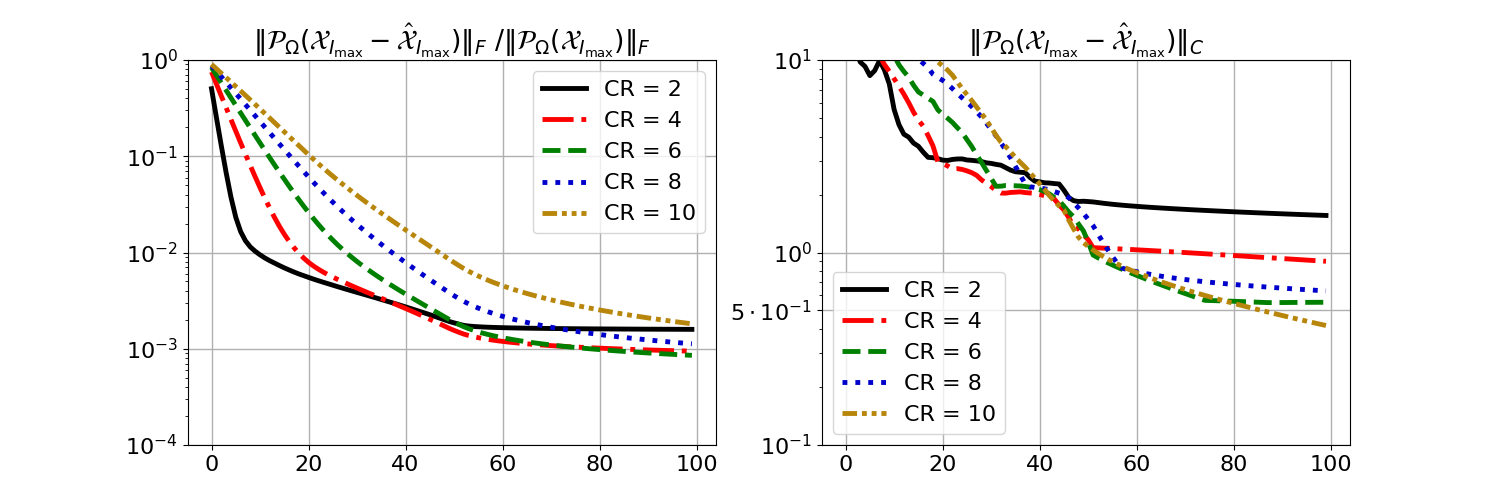}
\caption{Convergence of the SVP algorithm on the random mask of the measured elements, assessed by the Frobenius norm (left) and Chebyshev norm (right) for $\mathcal{X}_{I_{\max}}$ of May data. The dependence on the number of iterations is shown for different compression ratios $CR = 2, 4, \dots, 10$..}
    \label{fig:SVP_convergence_may}
\end{figure}

The right panel in Figures \ref{fig:SVP_convergence_jan} and \ref{fig:SVP_convergence_may} marks an acceptable error threshold of $|\mathcal{P}{\Omega}(\mathcal{X}{I_{\text{max}}}-\mathcal{\hat{X}}{I{\text{max}}})|=0.5$, providing a practical benchmark for assessing convergence speed.

This performance is contrasted with the May dataset results in Figure \ref{fig:SVP_convergence_may}. The convergence profiles are similar; however, both error norms are consistently higher across all compression ratios. We hypothesize that this degradation in performance is directly attributable to the increased dynamical complexity and faster rate of temperature change in the Black Sea during the spring season. This higher variability makes the data inherently less compressible via the completion techniques, presenting a greater challenge for the low-rank reconstruction model and resulting in a higher residual error upon convergence.

The convergence profiles from Figures \ref{fig:SVP_convergence_jan} and \ref{fig:SVP_convergence_may} demonstrate that the SVP-inspired successfully minimizes the error metrics given the constraints of the random mask $\mathcal{P}_{\Omega}(\cdot)$. However, the ultimate completion errors, quantified in Table \ref{tab:completion_errors_combined}, reveal that the overall quality of the compression is poor and unacceptable for practical applications.

For instance, even at a low compression ratio of CR=2, the relative Frobenius error for January is 0.0044 and the Chebyshev norm error exceeds 13 that is unacceptable for the modelling and processing purposes. These errors grow substantially with increasing CR. The results for May are markedly worse, with errors nearly double those of January across all CRs, consistent with the convergence analysis.

\begin{table}[H]
\begin{center}
\caption{Completion errors of the largest block for January and May data with respect to the different compression ratios. The error in Chebyshev norm is unacceptable even though the relative error in the Frobeius norm is quite small.}
\vspace{0.6cm}
\label{tab:completion_errors_combined}
\begin{tabular}{| c | c c | c c |}
\hline
\multirow{2}{0.8cm}{CR} & \multicolumn{2}{c|}{\textbf{January}} & \multicolumn{2}{c|}{\textbf{May}} \\
\cline{2-5}
 & $\frac{\|\mathcal{X} - \mathcal{\hat{X}}\|_{F}}{\|\mathcal{X}\|_{F}}$ & $\|\mathcal{X} - \mathcal{\hat{X}}\|_{C}$ & $\frac{\|\mathcal{X} - \mathcal{\hat{X}}\|_{F}}{\|\mathcal{X}\|_{F}}$ & $\|\mathcal{X} - \mathcal{\hat{X}}\|_{C}$ \\
\hline
2  & 0.0044 & 13.14 & 0.0088 & 44.14  \\
4  & 0.0177 & 23.70 & 0.0329 & 54.90  \\
6  & 0.0371 & 26.79 & 0.0625 & 49.84  \\
8  & 0.0617 & 31.68 & 0.0968 & 50.23 \\
10 & 0.0924 & 26.12 & 0.1312 & 48.54 \\
\hline
\end{tabular}
\end{center}
\end{table}

This poor performance leads to the conclusion that the completion frameworks, and by extension other similar tensor completion algorithms based on random observation masks, is unsuitable for the effective compression of this specific dataset. The fundamental issue lies in the data's structure. Although the initial tensor $\mathcal{X}$ is technically sparse (contains missing values), the existing data is organized in large, contiguous spatio-temporal blocks rather than being randomly scattered. Completion algorithms are designed to infer missing points from randomly distributed neighbors. They struggle significantly to reconstruct the continuous structures present in this environmental data.

Hence, we infer that for this type of data, deterministic compression methods (e.g., based on truncating the Tucker or TT- decompositions) would be a far more suitable and effective choice than completion-based approaches.

\subsection{Tucker experiment}

In this section we present the results of numerical experiments on the application of Tucker decomposition for approximation and compression of the four-dimensional array of reanalysis data on the temperature of the Black, Azov and Marmara Seas. 

For the initial experiment, we select the largest data block $\mathcal{X}_{I_{\text{max}}}$ (see in Figure~\ref{fig:largest_block}) to assess the compression capabilities of Tucker decomposition. We employ the HOSVD algorithm to compute approximations with controlled multilinear ranks.

\begin{figure}[H]
    \centering
    \includegraphics[width=0.7\linewidth]{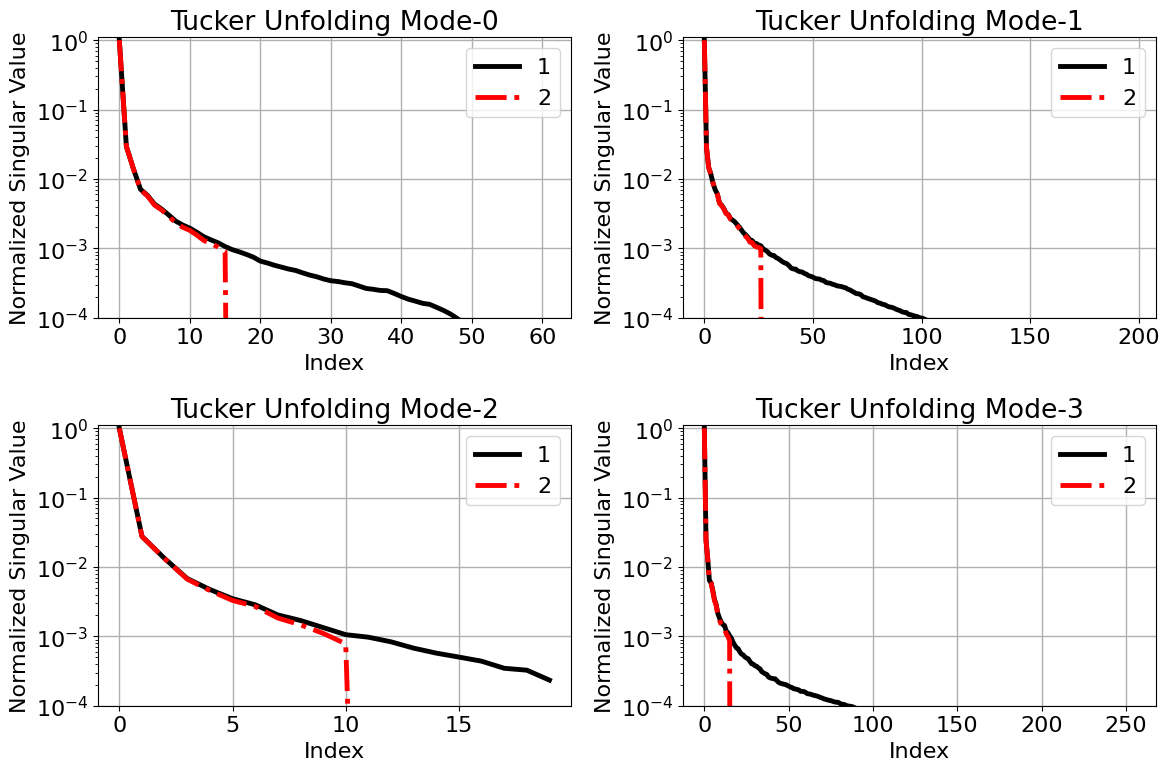}
    \caption{Normalized singular values of the four matrix unfoldings in Tucker decomposition. The singular values are normalized by division with the maximum singular value. Solid lines: complete spectra; dashed-dotted lines: values remaining after truncation with threshold $\varepsilon = 10^{-3}$.}
    \label{fig:HOSVD_truncation}
\end{figure}

Figure~\ref{fig:HOSVD_truncation} displays the normalized singular values for all four matrix unfoldings of $\mathcal{X}_{I_{\text{max}}}$, which has a temporal dimension of 256 time steps. We observe the rapid decay in these values indicating the strong compressibility potential through rank truncation.

Table~\ref{tab:compression_results_HOSVD_jan} shows the dependence of the compression ratio on the Tucker ranks and relative accuracy (Frobenius norm) for January and May respectively. The second column displays the Chebyshev error of the obtained approximation.

\begin{table}[H]
\centering
\caption{Tensor compression with HOSVD for different approximation rates. Experiment with the largest block of size $62\times 199\times 20\times 256$, January.}
\label{tab:compression_results_HOSVD_jan}
\vspace{0.6cm}
\begin{tabular}{c c c c c c c}
\hline
{$\varepsilon$} & {$\|\cdot\|_{\text{F}}$} & {$\|\cdot\|_{\text{C}}$} & {Ranks} & {Elements} & {Elements} & {Compression} \\
& & & & {before compression} & {after compression} & {rate (CR)} \\
\hline
$1\cdot 10^{-3}$ & 0.00523 & 2.2663 & [15, 25, 10, 15] & 63170560 & 66195 & 954.3 \\
 $5\cdot 10^{-4}$ & 0.00315 & 1.74746 & [23, 39, 15, 23] & 63170560 & 324840 & 194.5 \\
 $2\cdot 10^{-4}$ & 0.00151 & 0.96808 & [40, 70, 20, 44] & 63170560 & 2492074 & 25.3 \\
 $1\cdot 10^{-4}$ & 0.00079 & 0.72323 & [48, 94, 20, 80] & 63170560 & 7261762 & 8.7 \\
\hline
\end{tabular}
\end{table}

To study the dependence of the results on the season, a similar experiment for May was conducted. The results are presented in the table~\ref{tab:compression_results_HOSVD_May}. From the tables it can be seen that the data for January allow a more efficient tensor approximation than the data for May, which may be due to the greater variability of the Black Sea temperature in spring.

\begin{table}[H]
\centering
\caption{Tensor compression with HOSVD for different approximation rates. Experiment with the largest block of size $62\times 199\times 20\times 256$, May.}
\label{tab:compression_results_HOSVD_May}
\vspace{0.6cm}
\begin{tabular}{c c c c c c c}
\hline
{$\varepsilon$} & {$\|\cdot\|_{\text{F}}$} & {$\|\cdot\|_{\text{C}}$} & {Ranks} & {Elements} & {Elements} & {Compression} \\
 & & & & {before compression} & {after compression} & {rate (CR)} \\
\hline
$1\cdot 10^{-3}$ & 0.00807 & 3.85839 & [26, 36, 14, 22] & 63170560 & 302976 & 208.5 \\
 $5\cdot 10^{-4}$ & 0.00503 & 3.00241 & [43, 62, 18, 46] & 63170560 & 2234588 & 28.3 \\
 $2\cdot 10^{-4}$ & 0.00184 & 2.24747 & [52, 123, 20, 105] & 63170560 & 13486581 & 4.7 \\
 $1\cdot 10^{-4}$ & 0.00087 & 1.09985 & [59, 157, 20, 140] & 63170560 & 26007541 & 2.4 \\
\hline
\end{tabular}
\end{table}

Let us assume that the permissible absolute error in compressing temperature data is $0.5^{\circ}C$ \cite{satellites_description}. In the second numerical experiment, we compress all valid data blocks with an absolute error not exceeding $0.5^{\circ}C$. This condition is achieved through iterative selection of the appropriate multilinear ranks for the Tucker decomposition. For each block, we begin with conservative rank values and compute the Tucker approximation using the HOSVD algorithm. The resulting absolute error is then evaluated. If the absolute error exceeds the permissible threshold, the ranks are incrementally increased, and the approximation is recomputed. This iterative process continues until the absolute error falls within the acceptable range.

In Table~\ref{tab:detailed_compression_tucker} we present the approximation errors and compression ratios for the individual blocks, where the block numbers correspond to decreasing sizes (see Figure~\ref{fig:greedy_hmap}). This table includes the results for both January and May temperature data. The final row shows the overall compression ratios $CR_{all}$, calculated accounting for all data values including those outside the compressed blocks. Additionally, we report $CR_{\text{sub}}$, which measures compression efficiency relative to $\sum_{I_k} \mathcal{X}_{I_k}$ to isolate the effects of partitioning non-optimality from the Tucker decomposition performance itself.

\begin{table}[H]
\begin{center}
\caption{Results of Black and Azov Seas compression with Tucker and maximum absolute error $\varepsilon_{\max} = 0.5$}
\vspace{0.6cm}
\label{tab:detailed_compression_tucker}
\begin{tabular}{| c | c c c | c c c |}
\hline
\multirow{2}{1cm}{block no.} &  \multicolumn{3}{|c|}{\textbf{January}} & \multicolumn{3}{|c|}{\textbf{May}} \\
 & $\|\mathcal{X}_I-\hat{\mathcal{X}_I\|_{\text{F}}}/\|\mathcal{X}_I\|_{\text{F}}$ & ${\|\mathcal{X}_I-\hat{\mathcal{X}_I}\|_{\text{C}}}$ & ${CR_{i}}$ & $\|\mathcal{X}_I-\hat{\mathcal{X}_I\|_{\text{F}}}/\|\mathcal{X}_I\|_{\text{F}}$ & ${\|\mathcal{X}_I-\hat{\mathcal{X}_I}\|_{\text{C}}}$ & ${CR_{i}}$ \\
\hline
1 & 0.00047 & 0.46311 & 4.9 & 0.00035 & 0.48580 & 1.5 \\ 
2 & 0.00121 & 0.47869 & 31.2 & 0.00059 & 0.45812 & 2.1 \\
3 & 0.00063 & 0.45802 & 2.2 & 0.00113 & 0.46727 & 2.3 \\ 
4 & 0.00612 & 0.46922 & 20.8 & 0.00073 & 0.47070 & 4.5 \\ 
5 & 0.00073 & 0.49996 & 5.2 & 0.00108 & 0.44795 & 2.7 \\  
6 & 0.00139 & 0.47057 & 31.1 & 0.00093 & 0.48344 & 3.3 \\ 
7 & 0.00134 & 0.45865 & 4.8 & 0.00116 & 0.46024 & 2.3 \\ 
8 & 0.00062 & 0.45559 & 4.5 & 0.00150 & 0.49987 & 2.6 \\ 
9 & 0.00190 & 0.48313 & 5.4 & 0.00132 & 0.45891 & 2.4 \\
10 & 0.00217 & 0.48396 & 78.8 & 0.00138 & 0.40427 & 6.8 \\  
11 & 0.00410 & 0.47035 & 11.9 & 0.00149 & 0.46355 & 5.6 \\ 
12 & 0.00181 & 0.33624 & 15.8 & 0.00126 & 0.47409 & 2.3 \\ 
13 & 0.00094 & 0.47560 & 1.2 & 0.00120 & 0.46114 & 1.8 \\ 
14 & 0.00200 & 0.47512 & 5.1 & 0.00131 & 0.46377 & 2.6 \\ 
15 & 0.00380 & 0.45235 & 6.3 & 0.00217 & 0.49964 & 4.7 \\ 
16 & 0.01980 & 0.45429 & 187.6 & 0.00081 & 0.45394 & 7.2 \\ 
17 & 0.00183 & 0.47615 & 8.2 & 0.00235 & 0.45235 & 2.6 \\ 
18 & 0.00488 & 0.46854 & 8.9 & 0.00163 & 0.49160 & 6.0 \\ 
19 & 0.00359 & 0.47281 & 59.6 & 0.00220 & 0.48325 & 3.0 \\ 
20 & 0.04245 & 0.39328 & 125.6 & 0.00096 & 0.46878 & 15.3 \\ 
21 & 0.00155 & 0.44985 & 9.4 & 0.00150 & 0.47229 & 3.3 \\  
22 & 0.00282 & 0.39136 & 21.9 & 0.00179 & 0.42083 & 8.3 \\ 
23 & 0.00171 & 0.46061 & 10.9 & 0.00161 & 0.45299 & 5.4 \\  
\hline
~ & \multicolumn{6}{| c |}{\textbf{Total compression ratios:}} \\
~ & \multicolumn{3}{|c|}{CR\textsubscript{all} = 3.16} & \multicolumn{3}{|c|}{CR\textsubscript{all} = 1.67} \\
~ & \multicolumn{3}{|c|}{CR\textsubscript{cubes} = 5.37} & \multicolumn{3}{|c|}{CR\textsubscript{cubes} = 1.91} \\
\hline
\end{tabular}
\end{center}
\end{table}

This experiment investigates the impact of temporal partitioning on compression efficiency using Tucker decomposition. Preliminary analysis reveals that compression ratios are highly sensitive to the block size, particularly along the temporal dimension. For instance, even adjacent eight-day intervals within January may exhibit markedly different compressibility. In order to get a systematic evaluation of this dependency, we conduct experiments on a 32-day dataset compressed via the HOSVD procedure while varying the temporal division.

In Table~\ref{tab:cr_results_tucker_jan} (January) and Table~\ref{tab:cr_results_tucker_may} (May) we present the compression results under different partitioning schemes. The ``No Split'' configuration corresponds to the processing of the entire period as a single tensor, while ``N Splits" divides the temporal dimension into N equal intervals (with N being powers of two). Each resulting subtensor has to be compressed independently, with overall ratios calculated across the complete time domain. All experiments maintain the $0.5^{\circ}C$ maximal absolute error threshold.

\begin{table}[H]
\centering
\caption{Compression Ratio (CR) for different time splitting configurations (January)}
\vspace{0.6cm}
\label{tab:cr_results_tucker_jan}
\begin{tabular}{@{}lccccccc@{}}
\toprule
 & No Split & 2 Splits & 4 Splits & 8 Splits & 16 Splits & 32 Splits & 64 Splits\\
\midrule
$CR_{\text{all}}$ & 3.17 & 3.80 & 4.00  & 4.51 & 4.79 & 4.86 & 4.86 \\
$CR_{\text{sub}}$ & 5.37 & 8.10 & 9.28 & 13.39 & 16.89 & 18.05 & 18.05 \\
\bottomrule
\end{tabular}
\end{table}

\begin{table}[H]
\centering
\caption{Compression Ratio (CR) for different time splitting configurations (May)}
\vspace{0.6cm}
\label{tab:cr_results_tucker_may}
\begin{tabular}{@{}lccccccc@{}}
\toprule
 & No Split & 2 Splits & 4 Splits & 8 Splits & 16 Splits & 32 Splits & 64 Splits \\
\midrule
$CR_{\text{all}}$ & 1.67 & 1.62  & 1.64  & 1.76  & 1.87 & 2.01 & 2.06  \\
$CR_{\text{sub}}$ & 1.91 & 1.84  & 1.87  & 2.05  & 2.24 & 2.49 & 2.58 \\
\bottomrule
\end{tabular}
\end{table}

The results reveal a clear relationship between temporal partitioning granularity and compression efficiency. The compression ratios initially improve with finer partitioning, reaching an optimal value at two-day intervals. Beyond this point, the further subdivision leads to the performance degradation. This pattern suggests the existence of some optimal temporal scale for the Tucker decomposition of the ocean temperature data.

Consistent with the previous experiments, the May data exhibits significantly lower compression ratios compared to the January data, reinforcing the observed seasonal dependence.

\subsection{TT-experiment}

In this section we present the results of numerical experiments on the application of a tensor train (TT) for approximation and compression of our four-dimensional sea temperature data. In order to compare the performance of the Tucker and TT representations, we choose the largest block chosen in the first experiment (see Figure~\ref{fig:largest_block}). We use the TTSVD method \cite{oseledets2011tensor} to find the approximation with a given relative accuracy. 

Table~\ref{tab:compression_bb_TT_jan} shows the dependency of the compression ratio on the TT ranks and the relative accuracy (in Frobenius norm) for the January data. The second column shows the absolute error of the constructed approximation. These results can be compared with data from Table~\ref{tab:compression_results_HOSVD_jan}.

\begin{table}[ht]
\centering
\caption{Tensor compression with TT for different approximation rates. Experiment for TT with the largesrt block of size $62\times 199\times 20\times 256$, January}
\label{tab:compression_bb_TT_jan}
\vspace{0.6cm}
\begin{tabular}{ c c c c c c}
\hline
 {$\|\cdot\|_{\text{F}}$} & {$\|\cdot\|_{\text{C}}$} & {Ranks} & {Elements} & {Elements} & {Compression} \\
 & & & {before compression} & {after compression} & {rate (CR)} \\
\hline
 0.00458 & 2.214 & [16, 53, 13]   & 63\,170\,560 & 186\,852    & 338.1 \\
 0.00236 & 1.214 & [26, 143, 27]   & 63\,170\,560 & 825\,626   & 76.5  \\
 0.00119 & 0.751 & [37, 367, 60]   & 63\,170\,560 & 3\,160\,275   & 20.0   \\
 0.00058 & 0.515 & [45, 720, 105]  & 63\,170\,560 & 7\,989\,270  & 7.9   \\
 0.00029 & 0.327 & [50, 1125, 153]  & 63\,170\,560 & 14\,678\,518  & 4.3   \\
 0.00014 & 0.143 & [54, 1520, 200]  & 63\,170\,560 & 22\,468\,468  & 2.8   \\
 0.00007 & 0.090 & [57, 1869, 232]  & 63\,170\,560 & 29\,935\,153  & 2.1   \\
\hline
\end{tabular}
\end{table}
To study the dependence of the results on the season, a similar experiment for May was conducted. The results are presented in the Table~\ref{tab:compression_bb_TT_May}. From the tables it can be seen that the data for January allow a more efficient tensor approximation than the data for May, which may be due to the greater seasonal variability of the Black Sea temperature.
\begin{table}[ht]
\centering
\caption{Tensor compression with TT for different approximation rates. Experiment for TT with the largest block of size $62\times 199\times 20\times 256$, May}
\label{tab:compression_bb_TT_May}
\vspace{0.6cm}
\begin{tabular}{ c c c c c c}
\hline
 {$\|\cdot\|_{\text{F}}$} & {$\|\cdot\|_{\text{C}}$} & {Ranks} & {Elements} & {Elements} & {Compression} \\
 & & & {before compression} & {after compression} & {rate (CR)} \\
\hline
 0.00481 & 3.042 & [34, 224, 39]   & 59\,222\,400 & 1\,701\,772    & 34.8 \\
 0.00234 & 1.941 & [45, 516, 84]   & 59\,222\,400 & 5\,510\,610   & 10.8  \\
 0.00117 & 1.252 & [51, 849, 124]   & 59\,222\,400 & 10\,754\,943   & 5.5   \\
 0.00058 & 1.047 & [57, 1197, 167]  & 59\,222\,400 & 17\,619\,165  & 3.4   \\
 0.00027 & 0.484 & [61, 1522, 202]  & 59\,222\,400 & 24\,676\,700  & 2.4   \\
 0.00012 & 0.105 & [62, 1797, 223]  & 59\,222\,400 & 30\,243\,370  & 2.0   \\
 0.00006 & 0.073 & [62, 2012, 234]  & 59\,222\,400 & 34\,300\,220  & 1.7   \\
\hline
\end{tabular}
\end{table}

As in the Tucker experiment presented above, we set the acceptable absolute error in compressing the temperature data as $0.5^{\circ}C$ \cite{satellites_description}. In the second numerical experiment, we compress all valid data blocks with an absolute error not exceeding $0.5^{\circ}C$. This condition can be achieved by setting the required relative Frobenius approximation error quite small but still allowing to adapt the TT-ranks. Thus, for each block, a sufficiently large relative error is specified, and we construct the approximation using the TTSVD algorithm \cite{oseledets2011tensor}. Further, we calculate the resulting absolute error. If the absolute error exceeds the permissible level, the relative Frobenius error has to be reduced twice, and we recompute the approximation again. We repeat this process until the absolute error reaches the acceptable prescribed value. 

We show the approximation errors and the compression coefficients for each of the blocks (the first column shows the block number) in Table~\ref{tab:detailed_compression_tt}; the higher the number, the smaller the block size (see Figure~\ref{fig:greedy_hmap}). We demonstrate the compression results for the temperature in January and in May, and the last line in Table~\ref{tab:detailed_compression_tt} corresponds to the total compression ratios $CR_{all}$ obtained in the experiment. The values that do not belong to the blocks are also taken into account when calculating the coefficients $CR_{all}$ (these values are not compressed). 

Once again, we can see that the compression ratio of the different blocks differs significantly for the TT format as well as for the Tucker decomposition. We also see again that the compression rates for the temperatures in May are significantly lower than in January.

\begin{table}[ht!]
\begin{center}
\caption{Results of Black and Azov Seas compression with TT and maximum absolute error $\varepsilon_{\max} = 0.5$}
\vspace{0.6cm}
\label{tab:detailed_compression_tt}
\begin{tabular}{| c | c c c | c c c |}
\hline
\multirow{2}{1cm}{block no.} &  \multicolumn{3}{|c|}{\textbf{January}} & \multicolumn{3}{|c|}{\textbf{May}} \\
 & $\|\mathcal{X}_I-\hat{\mathcal{X}_I\|_{\text{F}}}/\|\mathcal{X}_I\|_{\text{F}}$ & ${\|\mathcal{X}_I-\hat{\mathcal{X}_I}\|_{\text{C}}}$ & ${CR_{i}}$ & $\|\mathcal{X}_I-\hat{\mathcal{X}_I\|_{\text{F}}}/\|\mathcal{X}_I\|_{\text{F}}$ & ${\|\mathcal{X}_I-\hat{\mathcal{X}_I}\|_{\text{C}}}$ & ${CR_{i}}$ \\
\hline
1 & 0.0004 & 0.428 & 5.73 & 0.0002 & 0.112 & 2.07 \\ 
2 & 0.0008 & 0.372 & 22.22 & 0.0004 & 0.444 & 4.20 \\
3 & 0.0009 & 0.381 & 3.65 & 0.0009 & 0.297 & 2.15 \\ 
4 & 0.0063 & 0.451 & 29.23 & 0.0008 & 0.425 & 2.40 \\ 
5 & 0.0008 & 0.445 & 10.30 & 0.0012 & 0.378 & 4.26 \\  
6 & 0.0008 & 0.359 & 17.51 & 0.0012 & 0.472 & 5.83 \\ 
7 & 0.0014 & 0.363 & 8.57 & 0.0012 & 0.411 & 3.22 \\ 
8 & 0.0007 & 0.441 & 8.20  & 0.0023 & 0.452 & 5.64 \\ 
9 & 0.0018 & 0.401 & 6.24 & 0.0010 & 0.283 & 1.85 \\
10 & 0.0015 & 0.381 & 49.35 & 0.0007 & 0.363 & 6.27 \\  
11 & 0.0026 & 0.303 & 6.22 & 0.0017 & 0.472 & 4.22 \\ 
12 & 0.0016 & 0.441 & 20.00 & 0.0010 & 0.296 & 1.88 \\ 
13 & 0.0006 & 0.399 & 1.09 & 0.0009 & 0.287 & 1.38 \\ 
14 & 0.0015 & 0.281 & 4.94 & 0.0012 & 0.308 & 2.94 \\ 
15 & 0.0037 & 0.362 & 4.27 & 0.0018 & 0.348 & 2.66 \\ 
16 & 0.0167 & 0.445 & 96.39 & 0.0008 & 0.255 & 2.32 \\ 
17 & 0.0014 & 0.276 & 7.92 & 0.0022 & 0.397 & 2.34 \\ 
18 & 0.0061 & 0.453 & 7.02  & 0.0016 & 0.376 & 3.25 \\ 
19 & 0.0028 & 0.368 & 37.82 & 0.0024 & 0.492 & 3.85 \\ 
20 & 0.0277 & 0.418 & 45.66 & 0.0003 & 0.180 & 3.10 \\ 
21 & 0.0013 & 0.282 & 11.31 & 0.0012 & 0.295 & 3.15 \\  
22 & 0.0025 & 0.339 & 12.00 & 0.0016 & 0.417 & 3.65 \\ 
23 & 0.0014 & 0.331 & 11.64 & 0.0012 & 0.244 & 4.49 \\  
\hline
~ & \multicolumn{6}{| c |}{\textbf{Total compression ratios:}} \\
~ & \multicolumn{3}{|c|}{CR\textsubscript{all} = 3.48} & \multicolumn{3}{|c|}{CR\textsubscript{all} = 1.997} \\
~ & \multicolumn{3}{|c|}{CR\textsubscript{cubes} = 11.10} & \multicolumn{3}{|c|}{CR\textsubscript{cubes} = 2.81} \\
\hline
\end{tabular}
\end{center}
\end{table}

\begin{figure}
\begin{center}
\begin{subfigure}{0.42\textwidth}
\includegraphics[width=0.8\textwidth]{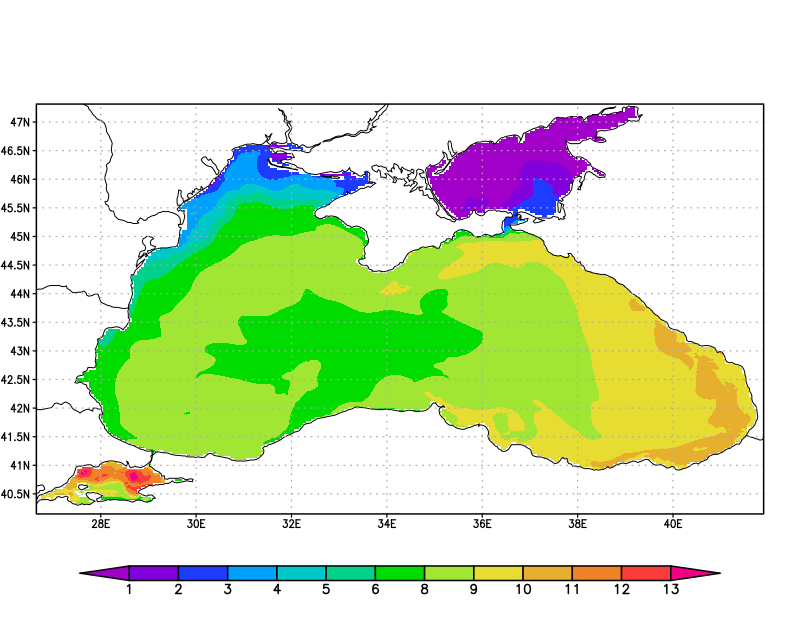}
\caption{~}
\end{subfigure}
\begin{subfigure}{0.42\textwidth}
\includegraphics[width=0.8\textwidth]{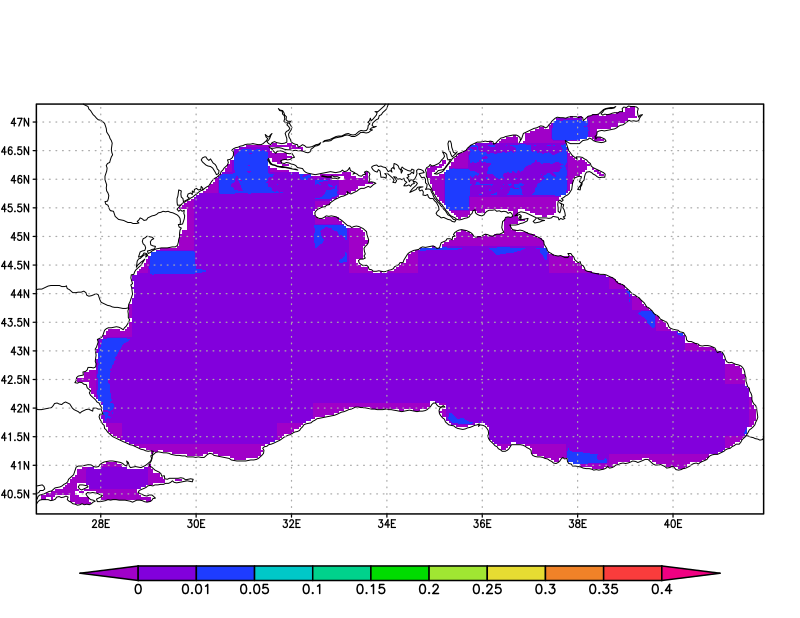}
\caption{~}
\end{subfigure}
\begin{subfigure}{0.42\textwidth}
\includegraphics[width=0.8\textwidth]{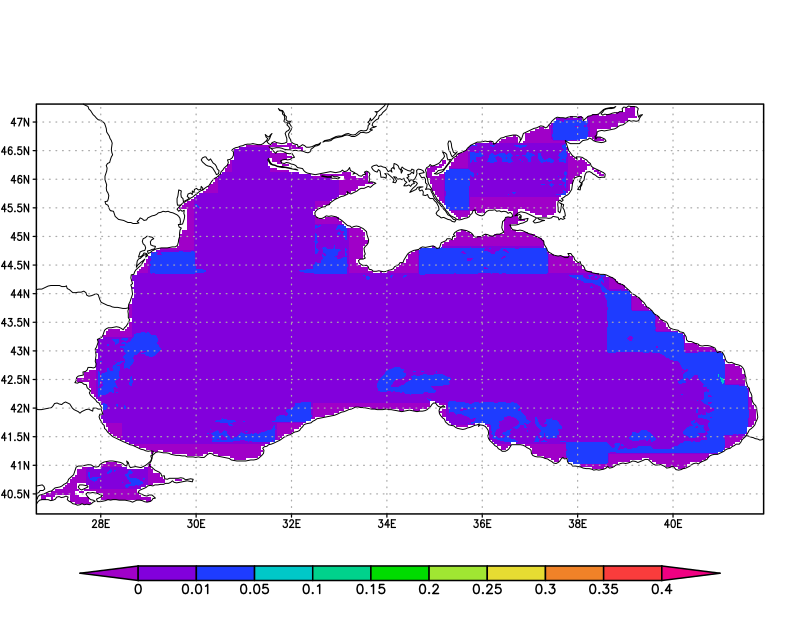}
\caption{~}
\end{subfigure}
\begin{subfigure}{0.42\textwidth}
\includegraphics[width=0.8\textwidth]{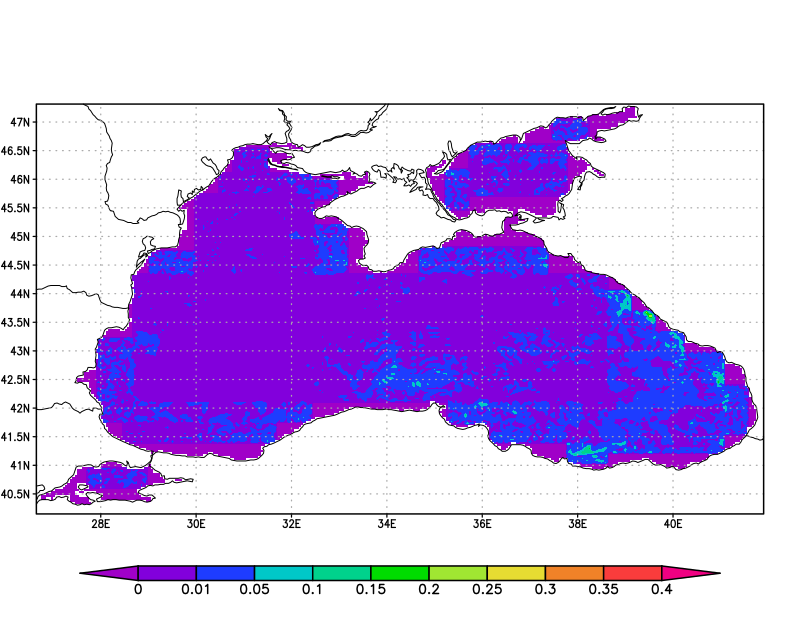}
\caption{~}
\end{subfigure}
    \caption{Numerical experiment on Black, Azov, and Marmara Seas temperature compression for January 2018 using TT. (a)  Sea surface temperature for January 22, restored after compression in TT format; (b) Localization of absolute approximation errors averaged over monthly time block and depth levels; (c) Localization of  absolute approximation errors averaged only over time for the surface level; (d) Absolute error of the sea surface temperature approximation for January 22. }
    \label{fig:errors_comparison}
\end{center}
\end{figure}

As we see from the experimental results presented in Tables \ref{tab:compression_bb_TT_jan}-\ref{tab:detailed_compression_tt}, where we are approximating the data tensor with an absolute error of 0.5 degrees Celsius, the relative Frobenius norm of the error remains quite small. However, it is also worth analyzing the spatial distribution of approximation errors, since the chosen compression approach may allow temperature ``jumps'' across block boundaries that may affect the quality of decompressed geophysical data. Figure \ref{fig:errors_comparison}(a) shows the sea surface temperature on January 22, reconstructed from the compressed data. Visually, the image is indistinguishable from the original data. Figure \ref{fig:errors_comparison}(d) shows the absolute difference between the original and compressed data for January 22. We can see from the figure that for this time layer the absolute approximation error does not exceed 0.2 degrees, with the errors located predominantly near the coast. It is also clear that the approximation error for the central (largest) block does not exceed 0.01 degrees for most points, that is a very good result. We select the layer corresponding to January 22 for demonstration, because it contains a fairly large number of points with errors, to show the characteristic structure of these errors. In Figure \ref{fig:errors_comparison}(c), the absolute error of the sea surface temperature approximation is averaged over all time layers in January, and in Figure \ref{fig:errors_comparison}(b), it is also averaged by depth.

In the next series of numerical experiments, we use different partitioning of data into the sub-tensors for comparison of our compression methods. Here we use the quantized tensor train (QTT) format \cite{oseledets2009approximation, sozykin2022ttopt} to compress the data. We construct the partitioning the same greedy algorithm as in the previous experiments (see Section 4), with the restriction that the sizes of the ``horizontal'' rectangles $x_{start}-x_{end}$, $y_{start}-y_{end}$ must be expressed as powers of two. We also use a different order of selecting the appropriate sizes here, one that prioritizes the same or similar sizes $x_{start}-x_{end}$, $y_{start}-y_{end}$. This allows to avoid choosing blocks of horizontal dimensions $4\times 128$, $64\times 2$, and so on. For the vertical dimension, prime factors are selected. For example, for $20$ of vertical levels in the considered dataset, the prime factors 2, 2 and 5 should be chosen. 

Let us consider as example, a data block $\mathcal X_I\in \mathbb R^{8 \times 8\times 20\times 8}$. We can approximate the original block in TT format or reshape the data and consider the multidimensional tensor ${\underbrace{2\times \dots \times 2}_{8}\times 5\times 2\times 2\times 2}$ (the so-called QTT procedure or introduction of the virtual dimensions into the data)\cite{oseledets2009approximation}. As the preliminary experiments have already shown, the data compression ratio is highly dependent on the block size. In particular, it is very sensitive to the choice of the time interval. For example, the data for the first eight days of January and the second eight days of January may be compressed differently. 

We also need to select and justify the optimal time interval for which data tensors are compressed most efficiently. In our experiment, the data for a period of 32 days are selected and compressed using TT and QTT. We present the results in Table~\ref{tab:cr_results_tt_jan} (the data for January) and Table~\ref{tab:cr_results_tt_may} (the data for May). The column ``No Split'' indicates that the entire 32-day dataset is compressed as a single tensor. Columns ``N Splits'', where $N$ is a power of two, indicate that the time dimension $T$ is divided into $N$ parts, and each $i$-th subtensor $\mathcal{X}\left(:,:,:,iT/N:(i+1)T/N\right)$ corresponding to i-th time subinterval, $i=0,1,\dots,N-1$, is compressed independently, while the compression ratio is calculated for the entire time interval. 

We summarize these observations with Fig. \ref{fig:errors_comparison}. Even though the compressed data in the TT-format is visually very close to the real sea temperature we see that the errors of approximation localize in different patterns closer to the coastlines of the sea. The averaging of these errors along the time and depth modes makes their pattern simpler supporting our motivation to apply the splitting and partitioning of the original data into the simple structured blocks.

\begin{figure}[H]
    \centering
    \includegraphics[width=0.49\linewidth]{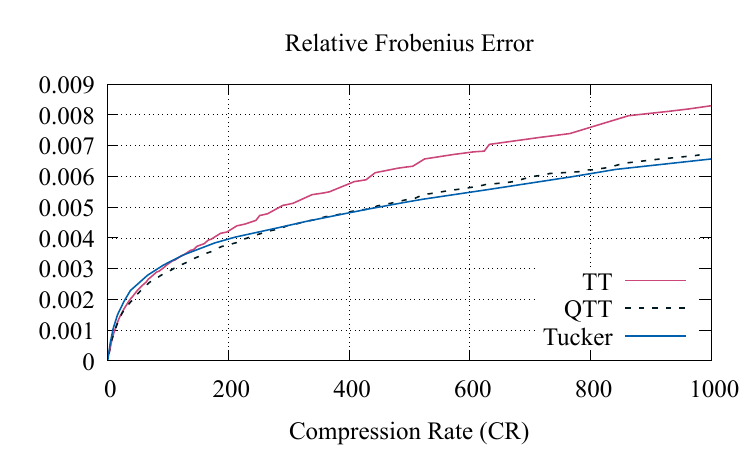}
    \includegraphics[width=0.49\linewidth]{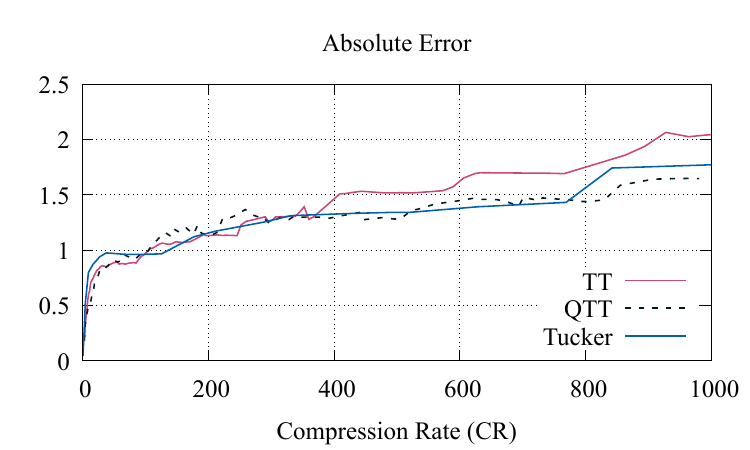}
    \caption{Left panel: Growth of the relative error in the Frobenius norm with respect to the growing compression rate. Right panel: Growth of the absolute error in the Chebyshev norm with respect to the growing compression rate. In contrast with application of the SVP-based completion the error in Chebyshev norm is does not exceed 2 degrees of Celsius even though the compression rate might be very large.}
    \label{fig:Compare_formats}
\end{figure}

From the results of the experiment we can see that the compression ratio is optimal when the data are divided into the subtensors corresponding to the two-day time intervals. The compression rates for the temperatures in May are still significantly lower than in January. This, utilization of the QTT-format does not provide significant benefits but allows to reach the same accuracy of approximation using the tensor train as Tucker decomposition (see Fig. \ref{fig:Compare_formats}) . In this experiment, the maximum permissible absolute error in compressing temperature data is chosen $0.5^{\circ}C$, as before.

\begin{table}
\centering
\caption{Compression Ratio (CR) for different time splitting configurations (January)}
\vspace{0.6cm}
\label{tab:cr_results_tt_jan}
\begin{tabular}{@{}lccccccc@{}}
\toprule
 & No Split & 2 Splits & 4 Splits & 8 Splits & 16 Splits & 32 Splits & 64 Splits\\
\midrule
$CR_{\text{all}}$ (TT) & 3.52 & 3.71 & 4.07  & 4.20 & 4.29 & 4.41 & 4.39 \\
$CR_{\text{all}}$ (QTT) & 3.34 & 3.83 & 4.02 & 4.20 & 4.38 & 4.41 & 4.30 \\
\bottomrule
\end{tabular}
\end{table}

\begin{table}[H]
\centering
\caption{Compression Ratio (CR) for different time splitting configurations (May)}
\vspace{0.6cm}
\label{tab:cr_results_tt_may}
\begin{tabular}{@{}lccccccc@{}}
\toprule
 & No Split & 2 Splits & 4 Splits & 8 Splits & 16 Splits & 32 Splits & 64 Splits\\
\midrule
$CR_{\text{all}}$ (TT) & 2.06 & 2.17 & 2.25  & 2.36 & 2.46 & 2.45 & 2.32 \\
$CR_{\text{all}}$ (QTT) & 2.07 & 2.19 & 2.31 & 2.49 & 2.59 & 2.51 & 2.34 \\
\bottomrule
\end{tabular}
\end{table}

\section{Conclusions}

This work set out to investigate efficient compression methods for the large tensors with gaps in their entries representing Black, Azov, and Marmara Seas temperature data. Initially, we develop a compression framework based on the Singular Value Projection (SVP) algorithm \cite{petrov2024constructing} operating under assumption that a tensor completion approach was necessary to handle the extensive formally missing values. However, in our experiments we obtain a surprising and clear result: the completion algorithm, while converging, produces the unacceptably high reconstruction errors, making it unsuitable for the practical application. The fundamental issue is a mismatch between the algorithm design for the randomly missing data and our dataset having a structure of large, contiguous gaps.

In contrary to initial expectations, we find a more direct and robust solution not in the complicated completion procedure, but in a strategic combination of data partitioning and standard SVD-based tensor decompositions. We introduce a greedy algorithm dividing the data into the spatially regular subtensors and create the basic units that are highly amenable to compression within the Tucker, TT, and the QTT formats. This approach effectively sidesteps the challenges posed by the original data gaps.

Our experiments confirm the effectiveness of this method. All three decomposition types (the Tucker, TT, and the QTT) achieve the high compression ratios while maintaining a strict maximum absolute error of 0.5°C, a feat the completion-based approach could not match. A key and unexpected finding is the strong temporal dependency of compressibility for our specific sea-temperature data. We observe the significant seasonal variation, with the January data consistently allowing for higher compression than the same data for May. Furthermore, we empirically determine that a temporal partitioning of approximately two days is nearly optimal for all decomposition methods, highlighting a kind of  universal temporal scale in the data structure.

In summary, the superior performance is achieved by a simpler but more efficient pipeline of intelligent data partitioning followed by standard SVD-based compression. This combined strategy proves to be a robust and effective method for managing this class of complex environmental data, achieving a consistent 4:1 compression ratio without compromising the accuracy required for scientific analysis.

\section*{Funding}

This work was partially supported (greedy partitioning method and experiments with tensor train format) by the Moscow Center of Fundamental and Applied Mathematics at INM RAS (Agreement with the Ministry of Education and Science of the Russian Federation No.075-15-2025-347); the algorithms and experiments with Tucker format were supported by Russian Science Foundation  25-11-00392 (the project webpage is available \url{https://rscf.ru/project/25-11-00392/}).

\section*{Declaration of interests}

The authors (S. Matveev, T. Sheloput and I. Kosolapov) declare that they have no known competing financial interests or personal relationships that could have appeared to influence the work reported in this paper. 

\bibliography{main}
\bibliographystyle{ieeetr}

\end{document}